\documentclass[]{amsart}
\usepackage[foot]{amsaddr}
\usepackage[a4paper,top=2.5cm,bottom=3cm,inner=3cm,outer=3cm]{geometry}

\usepackage[utf8]{inputenc}

\usepackage{amsmath,amsthm}
\usepackage{amsfonts,amssymb,mathrsfs}

\usepackage{xy}
\usepackage[bookmarks=false,breaklinks=true]{hyperref}
\usepackage{float}
\usepackage{graphicx}
\usepackage[percent]{overpic}
\usepackage{color}

\usepackage{epigraph}
\definecolor{myurlcolor}{rgb}{0,0,0.7}
\usepackage{hyperref}
\hypersetup{colorlinks,
linkcolor=myurlcolor,
citecolor=myurlcolor,
urlcolor=myurlcolor}
\usepackage[mathscr]{euscript}

\input xy
\xyoption{all}
\newcommand{\maps}{\colon}    %correct symbol for colon in f: X -> Y
\newcommand{\M}{{\mathcal M}}   % the Gromov-Hausdorff space
\newcommand{\T}{{\mathcal T}}  % the set of triangulations of the 2-sphere

\renewcommand{\H}{\mathrm{H}}  % Hausdorff metric
\newcommand{\GH}{\mathrm{GH}} % Gromov-Hausdorff metric

\newcommand{\R}{{\mathbb R}}  %real numbers

\newcommand{\define}[1]{{\bf \boldmath{#1}}}

\theoremstyle{definition}

        \newcommand{\be}{\begin{equation}}
        \newcommand{\ee}{\end{equation}}
        \newcommand{\ba}{\begin{eqnarray}}
        \newcommand{\ea}{\end{eqnarray}}
        \newcommand{\ban}{\begin{eqnarray*}}
        \newcommand{\ean}{\end{eqnarray*}}
        \newcommand{\barr}{\begin{array}}
        \newcommand{\earr}{\end{array}}

\begin{document}
\title{The Brownian Map}
\author[Baez]{John C.\ Baez} 
\address{Department of Mathematics, University of California, Riverside CA, 92521, USA}
\address{Centre for Quantum Technologies, National University of Singapore, 117543, Singapore}
\date{\today}
\maketitle

Nina Holden won the 2021 Maryam Mirzakhani New Frontiers Prize for her work on random surfaces and the mathematics of quantum gravity in 2 dimensions.   Here I will just explain one concept involved in her work: the `Brownian map'.  This is a fundamental object in mathematics, in some sense a 2-dimensional analogue of Brownian motion.  

Suppose you randomly choose a triangulation of the 2-sphere with $n$ vertices.  This is defined in a purely topological way, but the set of vertices becomes a metric space if we give each edge length 1.    We thus obtain a ``random compact metric space'': a probability measure on the space of all compact metric spaces.  It turns out that if you increase $n$ while rescaling the edge lengths to make them equal $n^{-1/4}$, this probability measure converges as $n \to \infty$.  The limit is a random compact metric space called the ``Brownian map''.      

To make all this precise takes a bit of work, and we turn to that next.  But first, here is a rough image to keep in mind, created by Thomas Budzinski:

\begin{center}
\includegraphics[width = 20em]{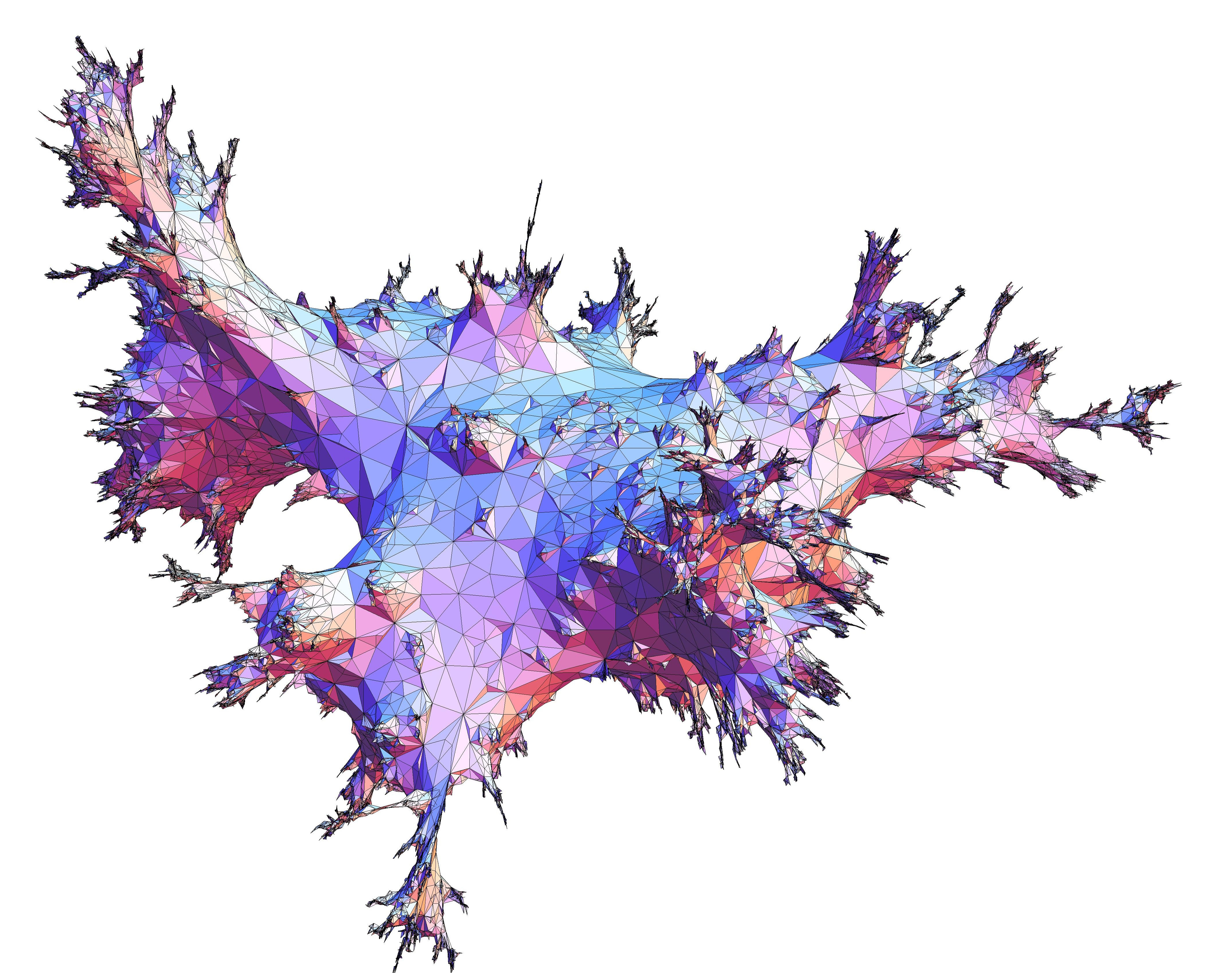}
\end{center}

\noindent
This shows a randomly chosen triangulation of the sphere with 30,000 vertices.  While this image gives an impression of the space's fractal structure, it is somewhat misleading, since the embedding in 3d space is not isometric.   This space is too spiky to be embedded isometrically in $\R^3$. Indeed, in 2007 Jean-Fran\c cois Le Gall \cite{LeGall07} showed that the Brownian map almost surely---that is, with probability 1---has Hausdorff dimension 4.  Nonetheless, it is almost surely homeomorphic to the 2-sphere \cite{LeGall08}.   For comparison, consider the famous Koch snowflake curve, which has Hausdorff dimension $\log_3(4) \approx 1.26$ but is homeomorphic to a 1-sphere.

Now let us make things more precise.  First, what do we mean by ``randomly choosing a triangulation of the 2-sphere with $n$ vertices?''   A triangulation of the 2-sphere is a way of gluing together a finite number of triangles along their edges, matching edges in pairs, so that the resulting space is homeomorphic to $S^2$.   The key issue is when two such triangulations count as the same.  The edges and vertices of any such triangulation give a graph topologically embedded in a surface homeomorphic to $S^2$, and we say two such triangulations count as the same if the resulting embedded graphs are related by an orientation-preserving homeomorphism between the surfaces containing them.    We obtain a finite set $\T_n$ of triangulations of the 2-sphere with $n$ vertices, and when we speak of choosing one randomly, we are putting a probability distribution on $\T_n$ where each element has the same probability.  Each element of $\T_n$ gives a compact metric space: namely, its set of vertices made into a metric space by taking all the edges of the graph to have length $n^{-1/4}$, and defining the distance between two vertices to be the length of the shortest path of edges between them.

Second, what is ``the space of all compact metric spaces''?   The \define{Hausdorff metric} on compact subsets of a metric space $K$ is defined by
\[    d_\H(X,Y) = \inf \{\epsilon \ge 0 : \; X \subseteq Y_\epsilon \textrm{\; and \;} Y \subseteq X_\epsilon \} \]
where $X_\epsilon$ is the set of all points of $K$ that are at distance at most $\epsilon$ from some point of $X$.  We define the \define{Gromov--Hausdorff distance} $d_{\GH}(X,Y)$ between compact metric spaces $X$ and $Y$ to be the infimum of $d_\H(i(X), j(Y))$ over all distance-preserving embeddings $i \maps X \to K$, $j \maps Y \to K$, where $K$ ranges over all compact metric spaces.   If there is an \define{isometry} between two compact metric spaces---that is, a distance-preserving bijection $f \maps X \to Y$---then clearly $d_{\GH}(X,Y) = 0$.  However, the Gromov--Hausdorff distance defines a metric on the set $\M$ of isometry classes of compact metric spaces.    The set $\M$ with its metric $d_{\GH}$ is called the \define{Gromov--Hausdorff space}.

Third, what is a \define{random compact metric space}?    It is simply a Borel measure $\mu$ 
on the Gromov--Hausdorff space that is a probability measure, meaning that $\int_\M \mu = 1$.  In particular, for each $n$, the probability distribution on the set $\T_n$, together with the way of making each triangulation of the 2-sphere into a compact metric space, gives a random compact metric space $\mu_n$.

Fourth and finally, what does it mean to say this sequence $\mu_n$ of random compact metric spaces converges to some random compact metric space $\mu$?    There are many concepts of convergence for probability measures, but here we want \define{convergence in distribution}: we require that for every bounded continuous function $f \maps \M \to \R$ we have
\[    \lim_{n \to \infty} \int_\M f d\mu_n \to \int_\M f d\mu  .\]

With these details filled in, it makes precise sense to ask if the random compact metric spaces $\mu_n$ converge as $n \to \infty$.    Oded Schramm \cite{Schramm} posed this as a problem in 2006, and Le Gall \cite{LeGall13} solved it in 2013: we have $\mu_n \to \mu$ for a unique random compact metric space $\mu$, the \define{Brownian map}.

But this result is just a small part of something much bigger: the slow march toward making ideas from quantum field theory mathematically rigorous.    Conformal field theory, a key tool in string theory, is the study of quantum fields on Riemann surfaces.  Any Riemann surface has a conformal structure, which lets us measure angles between tangent vectors, and a conformal field theory should transform in a simple way under automorphisms of Riemann surfaces---that is, conformal 
transformations.   The Brownian map arises naturally from a conformal field theory called ``Liouville quantum gravity''.    

To quickly get a feeling for this, suppose $\Sigma$ is a Riemann surface and $g$ a Riemannian metric compatible with the conformal structure on $\Sigma$.    Then there is a random distribution on $\Sigma$ called the \define{Gaussian free field}, defined by
\[  G = \sum_i g_i \phi_i  \]
where $g_i$ are independent Gaussian random variables with mean $0$ and variance $1$, and $\phi_i$ runs over an orthonormal basis of real eigenfunctions of the Laplacian on $\Sigma$ with nonzero eigenvalues---or more precisely, a maximal collection of such eigenfunctions with
\[ \frac{1}{2\pi} \int_\Sigma \nabla \phi_i(z) \cdot \nabla \phi_j(z) \, d^2 z  = \delta_{ij} .\]
Roughly speaking, Liouville quantum gravity is the study of the random Riemannian metric on 
$\Sigma$ given by
\[      g_\gamma =  e^{\gamma G} g  \]
where $\gamma \in \R$.   In reality things are not quite so simple.  For starters, $G$ is not a random function on $\Sigma$, merely a random distribution, so $g_\gamma$ is not a random smooth Riemannian metric.   Nonetheless, one can construct a random compact metric space based on Liouville quantum gravity---and when $\gamma = \sqrt{8/3}$, and $\Sigma = S^2$, this gives the Brownian map!   

This fact has been shown in a number of ways, starting with Jason Miller and Scott Sheffield \cite{MillerSheffield} around 2016.  More recently, Holden and Xin Sun showed how to go the other way, and recover the Gaussian free field on the sphere from large uniform triangulations \cite{HoldenSun}. 

The explanation of the number $\sqrt{8/3}$ would take us deeper into conformal field theory.   For other values of $\gamma$, we obtain a random metric space that has some other Hausdorff dimension $d_\gamma$ (almost surely).   No general formula for $d_\gamma$ is known, but intuitively we expect that $d_\gamma$ is an increasing function of  $\gamma$.  In 2018, this was proved by Jian Ding and Ewain Gwynne \cite{DG}.   

For more on Liouville quantum gravity and random surfaces, Gwynne's introduction is a good place to start \cite{Gwynne}.   Sheffield has a nice introduction to the Gaussian free field \cite{Sheffield}.   Finally, Gwynne, Holden and Sun have written a long review of more recent work \cite{GHS}, including that for which Holden won the Maryam Mirzakhani New Frontiers Prize.

\end{document}